\documentclass{amsart}
\usepackage{graphicx}
\usepackage{amssymb}
\newtheorem{thm}{Theorem}[section]
\newtheorem{cor}{Corollary}[section]
\newtheorem{lem}{Lemma}[section]
\newtheorem{prop}{Proposition}[section]
\theoremstyle{definition}
\newtheorem{defn}{Definition}[section]
\theoremstyle{remark}
\newtheorem{rem}{Remark}[section]
\numberwithin{equation}{section}
\begin{document}

\title{A construction of a skewaffine structure in Laguerre geometry}%
\author{Andrzej Matra\'{s}}%
\address{Andrzej Matra\'{s}\\
Department of Mathematics and Informatics\\
UWM Olsztyn\\
\.Zo\l nierska 14\\
10-561 Olsztyn\\
Poland}%
\email{matras@uwm.edu.pl}

%\thanks{}%

%\date{}%
%\dedicatory{}%
%\commby{}%
% ----------------------------------------------------------------
\begin{abstract}
J. Andre constructed a skewaffine structure as a group space of a
normally transitive group. In the paper this construction was used
to describe such an external structure associated with a point of
Laguerre plane. Necessary conditions for ensure that the external
structure is a skewaffine plane are given.
\end{abstract}
\maketitle
% ----------------------------------------------------------------

%%%%%%%%%%%%%%%%%%%%%%%%%%%%%%%%%%%%%%%%%%%%%%%%%%%%%%%%%%%%%%%%%%%%%%%%%%%%%%%%%%%%%%%%%%%%%%%%%%%%%%%%%
\section*{Introduction}
The internal incidence structure associated with a point $p$ of a
Benz plane consists of all points nonparallel to $p$ and, as
lines, all circles passing through $p$ (extended by all parallel
classes not passing through $p$). This is an affine plane, called
the derived affine plane in $p$. Natural question arises about a
characterization of the external structure in $p$
by some linear structure.\\
Among a wide class of noncommutative (in general) linear
structures constructed by J. Andre (cf. \cite{A}) there are
skewaffine planes which are good candidates to obtain the
characterization we are looking for. One of classical example of a
skewaffine plane is the set of circles of Euclidean plane with
centers of circles as basepoints (cf. \cite{A}). Under the weak
conditions Wilbrink (cf. \cite{W}) constructed in a fixed point of
Minkowski plane known as residual nearaffine plane - a skewaffine
plane such that any straight line intersects any nonparallel line
in exactly one point. In the examples (circles and hyperbolas) the
construction is based on the observation that any circle (resp.
hyperbola) has exactly one center which can be taken as a
basepoint of a line corresponding to a given conic. This center is
the image of our point $p$ in the inversion with respect to the
circle. In the case of Laguerre planes this construction cannot be
used, since the image of the point $p$ in the inversion is
parallel to $p$ and the inversion does not distinguished any point
which can be taken as a base point. It seems to be natural to
investigate a symmetry with two pointwiese fixed generators. Two
fixed generators do not define a symmetry of Laguerre plane. To
make the construction uniquely determinated we point out an
invariant pencil $\langle p,K\rangle$ of circles tangent in the
point $p$. A basepoint of a line corresponding to a circle (which
does not pass through $p$) we obtain as a point of tangency of the
circle with some circle of the pencil $\langle p,K\rangle$.
Opposite to M\"{o}bius and Minkowski planes basepoints belong to
corresponding circles and the residual skewaffine plane is not
determined by the point $p$ only, but by the pencil
$\langle p,K\rangle$.\\
A large class of regular skewaffine planes was given in \cite{A},
\cite{P} as the group space $\mathrm{V}(\mathbf{G})$ of a normally
transitive group $\mathbf{G}$. In such skewaffine planes lines are
obtained as a join of a basepoint and the orbit of an other point
with respect to the stabilizer of the
basepoint.\\
The starting point of our paper is a group $\mathbf{G}$ of
automorphisms of Laguerre plane such that $\mathrm{V}(\mathbf{G})$
contains a set of circle not passing through $p$ as lines. The
group $\mathbf{G}$ fixes points parallel to $p$ and the pencil
$\langle p,K\rangle$. Minimal conditions for the group
$\mathbf{G}$ are transitivity (cf. A1) and circular transitivity
for some circle of the pencil $\langle p,K\rangle$ (cf. A2). The
latter axiom A3 stands for an assumption that for any
$L\in\mathcal{C}$ with $p\notin L$ there exists exactly one
$M\in\langle p,K\rangle$ tangent to $L$. It guarantees that each
circle not passing through $p$ corresponds to some line of
$\mathrm{V}(\mathbf{G})$. We show that under the axioms A1, A2, A3
the group $\mathbf{G}$ contains of $\mathbb{L}$-translations and
$\mathbb{L}$-strains fixing the pencil $\langle p,K\rangle$. This
group is a subdirect product of the normal subgroup of
$\mathbb{L}$-translations and an arbitrary subgroup of
$\mathbb{L}$-strains with a fixed center (cf. Theorem \ref{t3.2}).
We point out an example of nonovoidal Laguerre plane which
satisfies our assumptions for some pencil $\langle p,K\rangle$
(cf. Remark \ref{r2.2}). Miquelian Laguerre planes of
characteristic distinct from 2 fulfils the axioms for any pencil (cf. Remark \ref{r2.1}).\\
We use skewaffine planes to get a characterization of the tangency
in Laguerre planes. We obtain the condition deciding whether a
circle through two points tangent to a circle can be constructed
(cf. Theorem \ref{t4.2}). As an application we also show that the
tangency points of circles of pencils $\langle p,K\rangle$ and
$\langle q,x\rangle$, where $q$ is parallel to $p$, form a circle
cf. Theorem \ref{t4.1}.
\section{Notations and basic definitions.}
A  \it Laguerre plane \/ \rm is a structure
$\mathbb{L}=(\mathcal{P},\mathcal{C},-)$, where
\begin{enumerate}
\item[---]
$\mathcal{P}$ is a set of {\it points} denoted by small Latin letters,
\item[---]
$\mathcal{C}\subset 2^{\mathcal{P}}$ is a set of
{\it circles} denoted by capital Latin letters,
\item[---]
$-$ is an equivalence  relation on ${\mathcal P}$
called {\it parallelity} (cf., for example, \cite{C}).
\end{enumerate}
We use the notation $-$ to avoid misunderstanding with the parallelity for
skewaffine planes.

The equivalence
classes of the relation $-$ we call \it generators \/ \rm and
denote by capital Latin letters. The following condition must be
satisfied:
\begin{enumerate}
\item[(1)] Three pairwise nonparallel points can be uniquely joined by a circle.
\item[(2)]  For every circle $K$ and any two nonparallel points $p\in K, q\notin K$ there is precisely one circle
$L$ passing through $q$ which touches $K$ in $p$ (i.e. $K\cap L=\{p\})$.
\item[(3)]  For any point $p$ and a circle $K$ there exists exactly one point q such that $p\parallel q$ and $q\in K$.
\item[(4)]  There exists a circle containing at least three points, but not all points.
\end{enumerate}
The unique generator containing a point $p$ we denote by
$\overline{p}$. If $a,b,c$ are points pairwise nonparallel the
unique circle containing points $a,b,c$ we denote by
$(a,b,c)^{\circ}$. The circle tangent to a circle $K$ and passing
through points $p,q$ ($p\in K, q\notin K,p-\!\!\!\!/ q$) we denote
by $(p,K,q)^{\circ}$. If $p\in K$ by the symbol $\langle
p,K\rangle$ we denote the pencil of circles tangent to $K$ in a
point $p$. If $x,y$ are nonparallel, then the set of circles
containing points $x,y$ we call \it pencil of circles with
vertexes \rm $x,y$ and denote by $\langle x,y\rangle$. For a point
$x$ and a circle $K$ the unique point of
$K$ parallel to $x$ (exists by the condition (3)) we denote by $xK$\\
The \it derived plane in a point \rm $p$ of Laguerre plane
$\mathbb{L}$ consist of all points not parallel to $p$ and, as
lines, all circles passing through $p$ (excluding $p$) and all generators not passing through $p$ (notation $\mathbb{A}_{p}$). This is an affine plane.\\
An automorphism of Laguerre plane is a permutation of the point
set which maps circles to circles (and generators to
generators). The automorphism $\phi$ is called \it central \rm if
there exists a fixed point $p$ such that $\phi$ induces a central
collineation of $\overline{\mathbb{A}_{p}}$ - the projective
extension of the derived affine plane
$\mathbb{A}_{p}$.\\
\it$\mathbb{L}$-translation \rm  is a central automorphism of Laguerre plane $\mathbb{L}$ which fixes points of
$\overline{p}$ and induces a translation of $\mathbb{A}_{p}$ for some point $p$. The group of translations which fixes
circles of a pencil $\langle p,K\rangle$ (respectively the familly of generators $\mathcal{G}$) is denoted by
$\mathbf{T}(\overline{p},K)$ (resp. $\mathbf{T}(\overline{p},\mathcal{G})$) (cf. \cite{H}).\\
\it$\mathbb{L}$-strain \rm with respect to a generator is a central automorphism of $\mathbb{L}$ which fixes points of
some generator $\overline{p}$ and circles of some pencil $\langle q,M\rangle$ ($p-\!\!\!\!\!/q$). The group of all
$\mathbb{L}$-strains fixing the points of $\overline{p}$ and the circles of the pencil $\langle q,M\rangle$ we denote
by $\Delta(\overline{p},q,M)$.  The involutory automorphism
which fixes pointwise generators $X,Y$ and a circle $M$ (not pointwise) we call \it Laguerre symmetry \rm and denote
by the
symbol $\mathrm{S}_{X,Y;M}$ ($\mathrm{S}_{X,Y;M}\in\Delta(\overline{p},q,M)$ for $p\in X$, $q\in Y$).\\
A group of central automorphisms is called {\it circular
transitive} if the extension to $\overline{\mathbb{A}_{p}}$ of its
restriction to $\mathbb{A}_{p}$ is linear transitive i.e. is
transitive on any line passing through the center. The group of
automorphisms of $\mathbb{L}$ is called $\langle
p,K\rangle$-{\it transitive} (resp. $\overline{p}$-transitive) if it
contains a circular transitive group $\mathbf{T}(\overline{p},K)$
(resp. $\mathbf{T}(\overline{p},\mathcal{G}$), cf. \cite{K}). The
group of automorphisms of $\mathbb{L}$ is called
$(\overline{p},\langle q,M\rangle)$-{\it transitive} if it contains a
circular
transitive group $\Delta(\overline{p},q,M)$.\\
A \it skewaffine space \rm (cf. \cite{A}) is an incidence
structure $\mathbb{S}=(X,\sqcup,\parallel)$, where $X$ is a
nonempty set of points denoted by small Latin letters, and
$$\sqcup:\{(x,y)\in X^{2}|x\neq y\}\rightarrow 2^{X}$$
is a function.
The sets of the form $x\sqcup y$ ($x\neq y$) are called lines.
They will be also denoted by capital Latin letters. The symbol $\parallel$
denotes an
equivalence relation among the lines. The following axioms must be satisfied:
\begin{enumerate}
\item[(L1)] $x,y\in x\sqcup y$,
\item[(L2)] $z\in x\sqcup y\setminus\{x\}$, implies $x\sqcup y=x\sqcup z$ (exchange condition),
\item[(P1)] given any line $L$ and any point $x$ there exists exactly one line $x\sqcup y$ parallel to $L$
(Euclidean axiom),
\item[(P2)] $\forall x,x',y,y'(x\neq y, x'\neq y'\wedge x\sqcup y\parallel x'\sqcup y'\rightarrow y\sqcup x\parallel y'\sqcup
x'$ (symmetry condition),
\item[(T)] if $x,y,z$ are pairwiese different points such that $x\sqcup y\parallel x'\sqcup y'$
there exists a points $z'$
such that $x\sqcup z\parallel x'\sqcup z'$ and  $y\sqcup z\parallel y'\sqcup z'$ (Tamaschke's condition).
\end{enumerate}
If we assume $x=x'$ in axiom (T) the condition is called affine
Veblen-condition (V). We will consider additional
conditions for skewaffine space:
\begin{enumerate}
\item[(Pgm)]
$\forall x,y,z\in X,\{x,y,x\}_{\neq}\exists w\in X$ with $x\sqcup y\parallel z\sqcup w$ and $x\sqcup z\parallel
y\sqcup w$
\item[(Des)]
$\forall u,x,y,z\in X,\{u,x,y,z\}_{\neq}x'\in u\sqcup x\setminus\{u\}\rightarrow\exists y'\in u\sqcup
y\setminus\{u\},z'\in u\sqcup z\setminus\{u\}$ with $x\sqcup y\parallel x'\sqcup y', x\sqcup z\parallel x'\sqcup
z',y\sqcup z\parallel y'\sqcup z'$
\item[(Pap)]
$\forall u,x,y,z,x'\in X$, $\{u,x,x'\}_{\neq}$ with $u\sqcup x=u\sqcup y=u\sqcup z\rightarrow\exists y',z'\in X$
with $u\sqcup x'=u\sqcup y'=u\sqcup z'$ and $x\sqcup x'\parallel z\sqcup z',x\sqcup y'\parallel y\sqcup z',y\sqcup
x'\parallel z\sqcup y'$.
\end{enumerate}
If a line $L$ takes the form $x\sqcup y$ then the point $x$ is called the basepoint of the line $L$. It is a simple consequence
of the axioms that any line has either exactly one basepoint or all its points are basepoints
 (cf. \cite{P}). A line all of its points are basepoints is called straight line. A line which is not straight (and hence has exactly
 one basepoint) is called a proper line.\\
 A group $\mathbf{G}$ acting on a set $X$ is called normally transitive if $\mathbf{G}$ is transitive and
 $\mathbf{G}_{x}\setminus\mathbf{G}_{y}\neq\emptyset$ hold for any $x,y\in X$ with $x\neq y$ ($\mathbf{G}_{x}$ denotes
 the stabilizer of a point $x$ with respect to $\mathbf{G}$). For any group acting on a set $X$ one can construct a
 group space $\mathbf{V}(\mathbf{G})=(X,\sqcup,\parallel)$ with\\
 $x\sqcup y =\mathbf{G}_{x}\{x,y\}=\{x\}\cup\mathbf{G}_{x}y$ and\\
 $L\parallel L'$ if there exists $g\in\mathbf{G}$ such that $gL=L'$ for any lines $L, L'$.\\
 The following theorem will be basic in our construction.
 \begin{thm}
 The group space $\mathbf{V}(\mathbf{G})$ with respect to a normally transitive group $\mathbf{G}$ is a desarguesian
 skewaffine space.
 \end{thm}
 More detailed discussion of properties of the group space $\mathbf{V}(\mathbf{G})$ can be found in \cite{P}.

\section{Residual skewaffine plane}
Let $\langle p,K\rangle$ be a fixed pencil of circles of Laguerre plane $\mathbb{L}=(\mathcal{P},\mathcal{C},-)$ and
$\Delta(p,K)$ a group of automorphisms of $\mathbb{L}$ which fixes the pencil $\langle p,K\rangle$ (not pointwiese) and
all the points parallel to $p$. We assume that $\Delta(p,K)$ satisfies the following conditions:
\begin{enumerate}
\item[(A1)] $\Delta(p,K)$ is transitive on the set $\mathcal{P}\setminus p$.
\item[(A2)] $\Delta(p,K)_{r}$ is circular transitive for any $r\in K\setminus\{p\}$ (i.e. for any $x,y\in K\setminus\{p,r\}$
there exists $\sigma\in\Delta(p,K)_{r}$ such that $\sigma(x)=y)$).
\end{enumerate}

\begin{defn}\label{d2.1}
A \it residual skewaffine plane with respect to \rm $\langle
p,K\rangle$ (written $\mathbb{SA}(p,K)$) is the group space
$\mathbf{V}(\Delta(p,K)=(\mathcal{P}\setminus
\overline{p},\sqcup,\parallel)$.
\end{defn}
By the definition $\Delta(p,K)$ is normally transitive and according to Theorem 1(\cite{A}, p. 5) (comp. also \cite{P}
Proposition 6.5, p. 94) we get
\begin{thm}\label{t2.1}
Suppose $\langle p,K\rangle$ is a fixed pencil of a Laguerre plane
$\mathbb{L}=(\mathcal{P},\mathcal{C},-)$ and conditions
$\mathrm{(A1), (A2)}$ are satisfied. Then the residual skewaffine
plane $\mathbb{SA} (p,K)$ is a skewaffine desarguesian space.
\end{thm}
In the following let $\langle p,K\rangle$ be a fixed pencil such that the
corresponding group $\Delta(p,K)$ satisfies the conditions (A1), (A2).
\begin{prop}\label{p2.1}
a) If $r-\!\!\!\!\!/x$, then $r\sqcup x=M\setminus \overline{p}$, where $M$ is a circle of $\langle p,K\rangle$ or a
circle tangent in the point $r$ to some circle of $\langle p,K\rangle$.\\
b) If $r-x$, then $r\sqcup x\subseteq \overline{r}$.
\end{prop}
\begin{proof}
a) If $r\notin K$, then consider a circle $L=(p,K,r)^{\circ}$.
First assume $x\in L$. If $r'$ is an arbitrary point of $K$
different from $p$ there exists an automorphism
$\alpha\in\Delta(p,K)$ such that $\alpha(r)=r'$. From the
definition of the group $\Delta(p,K)$, by the touching axiom we
get $\alpha(L)=K$. Hence, by (A2),
$K\setminus\{p\}=\Delta(p,K)_{r'}\alpha(x)$ and
$L\setminus\{p\}=\alpha^{-1}(K\setminus\{p\})=\Delta(p,K)_{r}x$.
In the case $x\notin L$ we define $M=(r,L,x)^{\circ}$, $x'=xL$ and
$y\in M\setminus\overline{p}$, $y\neq r$. Similarly as above we
conclude that $M$ is invariant with respect to the group
$\Delta(p,K)_{r}$. By the proved part of the proposition there
exists $\beta\in\Delta(p,K)_{r}$ such that
$\beta(x')=yL$. Hence $\beta(x)=y$.\\
b) It follows directly from the definition.
\end{proof}
\begin{defn}\label{d2.2}
The line $x\sqcup y$ of $\mathbb{SA}(p,K)$ is called \it special \rm if $x-y$.
\end{defn}
From the proof of Proposition \ref{p2.1} we get the generalization of the property (A2).
\begin{cor}\label{c2.1}
Let $r\notin\overline{p}$ and $M$ is invariant with respect to $\Delta(p,K)_{r}$. Then for any $x,y\in M$ with $x,y\neq
pM,r$ there exists $\sigma\in\Delta(p,K)_{r}$ such that $\sigma(x)=y$.
\end{cor}
\begin{prop}\label{p2.2}
The lines determined by the circles of the pencil $\langle p,K\rangle$ are straight lines.
\end{prop}
\begin{proof}
If $x,y$ are distinct points of a circle $L\in\langle p,K\rangle$, then $\Delta(p,K)_{x}(y)=L\setminus\{p\}$ by
Corollary \ref{c2.1}.
\end{proof}
The group $\Delta(p,K)$ fixes the generator $\overline{p}$
pointwiese. Hence by the definition of parallelity we get:
\begin{prop}\label{p2.3}
For any $M,N\in \mathcal{C}$ if $M\setminus\overline{p}\parallel L\setminus\overline{p}$, then $pM=pL$.
\end{prop}
Propositions \ref{p2.1}, \ref{p2.2}, \ref{p2.3} provide the
representation of basic notions of residual skewaffine plane of
Laguerre plane. However this representation is not complete. The
construction does not assume that any circle of the Laguerre plane
corresponds to some line of skewaffine plane.
Similarly, the statements inverse to
Propositions \ref{p2.2} and \ref{p2.3} do not hold
(as an example we can take any miquelian Laguerre plane of the characteristic 2).
To obtain this
representation more complete the following assumption will be needed throughout the further of the paper.\\
(A3) For any circle $M$ such that $p\notin M$ there exists exactly one circle $L\in \langle p,K\rangle$ tangent to $M$.
\begin{prop}\label{p2.4}
If $M$ is a circle such that $p\notin M$ then
$M\setminus\overline{p}=x\sqcup y$, where $x$ is the point of
tangency $M$ with the unique circle of the pencil
$\langle p,K\rangle$ and $y$ is any point of the circle $M$ different from $x$ and $pM$.
\end{prop}
\begin{proof}
According to (A3) there exists exactly one circle $L\in\langle p,K\rangle$ tangent to $M$. Let $x$ be its point of
tangency with $M$. The circle $L$ is invariant with respect to $\Delta(p,K)_{x}$, since $L\in\langle p,K\rangle$ and
$x\neq p$ is fixed. Hence $M$, as a circle tangent to an invariant circle and containing a fixed point $pM$ is invariant
with respect to $\Delta(p,K)_{x}$. The assertion follows from Corollary \ref{c2.1}.
\end{proof}
According to Proposition \ref{p2.4} in the case $x\sqcup
y=M\setminus\overline{p}$ the base point $x$ of the line $x\sqcup
y$ we will also call the \it base point \rm of the circle $M$.
\begin{prop}\label{p2.5}
If $M\setminus\overline{p}$ is a straight line, then $M\in\langle p,K\rangle$.
\end{prop}
\begin{proof}
Assume that $M\setminus\overline{p}=x\sqcup y=y\sqcup x$ for some distinct points $x,y\in M\setminus\overline{p}$ and
$M\notin\langle p,K\rangle$. Then $(x,M,p)^{\circ}$ and $(y,M,p)^{\circ}$ would be two distinct circles of the pencil
$\langle p,K\rangle$ tangent to $M$, a contradiction with (A3).
\end{proof}
\begin{prop}\label{p2.6}
If $pM=pL$ with $p\notin M,L$, then $M\setminus\overline{p}\parallel L\setminus\overline{p}$.
\end{prop}
\begin{proof}
From the Proposition \ref{p2.4} $M\setminus\overline{p}=x\sqcup y$ where $x$ is a point of tangency $M$ with some
$M'\in\langle p,K\rangle$ and $L\setminus\overline{p}=z\sqcup t$ where $z$ is a point of tangency $L$ with some
$L'\in\langle p,K\rangle$. By (A1) there exists $\sigma\in\Delta(p,K)$ such that $\sigma(x)=z$. We obtain
$\sigma(M')=L'$ and hence $\sigma(M)=L$ by the touching axiom.
\end{proof}
\begin{rem}\label{r2.1}
In the case of miquelian planes of the characteristic different from 2
any pencil $\langle p,K\rangle$ satisfies (A3) and the group
$\Delta(p,K)$ has properties (A1) and (A2). Additionally,
$\mathbf{V}(\Delta(p,K))$  fulfills (Papp) and (Pgm).
\end{rem}
\begin{rem}\label{r2.2}
Examples of nonmiquelian (even nonovoidal) planes satisfying (A1), (A2), (A3) one may obtain if we put in the
construction from \cite{AG} $f(x)=|x|^{r}$ ($r>1$) and $p=(\infty,0)$, $K=\{(x,0)|x\in K\}\cup\{(\infty,0)\}$.
Transformations of the form $x'=kx$; $y'=|k|^{r}y$ (identity for points $(\infty,a)$) form the stabilizer of a point
$(0,0)$.
\end{rem}
\begin{rem}\label{r2.3}
Axiom (A3) is satisfied for any pencil $\langle p,K\rangle$ of a
topological Laguerre plane of dimension 2 and 4 (it is a special
case of the solution of Apolonius problem for such planes cf.
\cite{S}).
\end{rem}
\section{A characterization of the group $\Delta(p,K)$.}
\begin{thm}\label{t3.1}
For any point $r\notin\overline{p}$ the stabilizer $\Delta(p,K)_{r}$ is a
$(\overline{p},\langle r,L\rangle)$-transitive
group of $\mathbb{L}$-strains $\Delta(\overline{p},r,L)$ containing $\mathrm{S}_{\overline{r},\overline{p};K}$,
where $L=(p,K,r)^{\circ}$.
\end{thm}
\begin{proof}
For any $\phi\in\Delta(p,K)_{r}$, $\phi(L)=L$, since $\Delta(p,K)$ fixes $\langle p,K\rangle$. Hence $\phi(M)=M$ for
all $M\in\langle r,L\rangle$, because $\overline{p}$ is poinwiese fixed. The $(\overline{p},\langle
r,L\rangle)$-transitivity of the group $\Delta(p,K)$ follows from Corollary \ref{c2.1}. To obtain
$\mathrm{S}_{\overline{r},\overline{p};K}$ let us consider an arbitrary $x\notin L\cup\overline{r}\cup\overline{p}$ and
circles $M=(r,L,x)^{\circ}$, $N=(p,L,x)^{\circ}$. According to (A3) circles $M$ and $N$ are not tangent, so there
exists a point $y$ with $y\neq x$, $y\in M\cap N$. By $(\overline{p},\langle r,L\rangle)$-transitivity there exist
$\psi\in\Delta(p,K)_{r}$ such that $\psi(x)=y$. Hence $\psi(N)=N$ and $\psi(y)=x$. It shows that
$\psi=\mathrm{S}_{\overline{r},\overline{p};K}$.
\end{proof}
\begin{prop}\label{p3.1}
The group $\Delta(p,K)$ is $\langle p,K\rangle$-transitive.
\end{prop}
\begin{proof}
Let $x,y\in K$ be such that $\{x,y,p\}_{\neq}$. Consider an arbitrary circle $M$ such that $M\cap K=\{x,y\}$ and the
base point $r$ of $M$. According to Theorem \ref{t3.1} there exist the symmetries
$\mathrm{S}_{\overline{r},\overline{p};K}$ and $\mathrm{S}_{\overline{x},\overline{p};K}$ and the superposition
$\mathrm{S}_{\overline{r},\overline{p};K}\circ\mathrm{S}_{\overline{x},\overline{p};K}$ is a translation which
transforms $x$ on $y$.
\end{proof}
From the proof of Proposition \ref{p3.1} we obtain:
\begin{cor} \label {c3.1}
For any distinct points $x,y$ such that $x,y\neq p$, $x,y\in R\in\langle p,K\rangle$ there exists $r\in R$ such that
$\mathrm{S}_{\overline{r},\overline{p};K}(x)=y$
\end{cor}
\begin{lem}\label{l3.1}
Any fixpoint free (outside $\overline{p}$) automorphism from $\Delta(p,K)$ is a translation.
\end{lem}
\begin{proof}
Suppose there exists $x$ such that $x-\!\!\!\!\!/\phi(x)$ for some
$\phi$ fulfilling the assumptions of the Lemma. Then
$\phi(x)-\!\!\!\!\!/\phi^{2}(x)$. If $x\neq\phi^{2}(x)$, then the
automorphism $\phi$ fixes the circle
$M=(x,\phi(x),\phi^{2}(x))^{\circ}$. If $x=\phi^{2}(x)$, then
$\phi$ fixes any circle $M\in\langle x,\phi(x)$. Hence $p\in M$
unless otherwise a base point of the circle $M$ is fixed. It means
that $\phi$ is a translation with an invariant pencil $\langle
p,M\rangle$. In the case $x-\phi(x)$ the assumption that there
exists $y$ such $y-\!\!\!\!\!/\phi(y)$ follows that the circle
$N=(y,\phi(y),\phi^{2}(y))^{\circ}$ is an invariant and a point of
intersection $N$ with $\overline{x}$ is fixed, a contradiction.
Thus in this case $\phi$ is a translation which fixes all
generators.
\end{proof}
\begin{prop} \label{p3.2}
The group $\Delta(p,K)$ is $\overline{p}$-transitive.
\end{prop}
\begin{proof}
Let $x-y$ and $r-\!\!\!\!\!/x$. Suppose $\phi\in\Delta(p,K)_{r}$ an arbitrary, $z=\phi(x)$,
$z'=\mathrm{S}_{\overline{p},\overline{r};K}(z)$ and $M=(y,z,z')^{\circ}$. Because
$\mathrm{S}_{\overline{p},\overline{r};K}(M)=M$ the base point $s$ of the circle $M$ is parallel to $r$. By Theorem
\ref{t3.1} there exists $\psi\in\Delta(p,K)_{s}$ such that $\psi(z)=y$, so $\psi\circ\phi(x)=y$. The automorphism
$\psi\circ\phi$ fixes (not pointwiese) two generators distinct from $\overline{p}$. Besides it is not an
$\mathbb{L}$-strain so it is fixpoint free outside $\overline{p}$. By Lemma \ref{l3.1} $\psi\circ\phi$ is a translation
which maps $x$ on $y$.
\end{proof}

The group of translations contained in $\Delta(p,K)$ will be denoted by $\mathbf{T}(p,\Delta)$.

\begin{thm}\label{t3.2}
Elements of the group $\Delta(p,K)$ without fixed points (outside
$\overline{p}$) are translations in direction of any $L$ $(p\in
L$) and translations which fix generators. The group $\mathbf{T}(p,\Delta)$
is transitive on the set $\mathcal{P}\setminus\overline{p}$ and
$\mathbf{T}(p,\Delta)\unlhd\Delta(p,K)$. Elements of $\Delta(p,K)$
with fixed points are $\mathbb{L}$-strains.
The group $\Delta(p,K)$ is isomorphic to $\simeq\mathbf{T}(p,\Delta)\rtimes \Delta(p,K)_{r}$ for
any $r\notin\overline{p}$.
\end{thm}
\begin{proof}
The first part of the theorem follows from Proposition \ref{p3.1}, Proposition \ref{p3.2}, and  \cite{HP}, Theorem
4.19. For any $\mathbb{L}$-strain $\phi$ and a translation $\psi$ the superposition $\phi\circ\psi\circ\phi^{-1}$ is a
translation. Indeed, otherwise $x=\phi\circ\psi\circ\phi^{-1}(x)$ for some $x\in\overline{p}$ and $\phi^{-1}(x)$ is the
fixed point of the translation $\psi$, a contradiction.
\end{proof}
\begin{cor}\label{c3.2}
The group $\Delta(p,K)$ is of type $\mathbf{1H}$ in the classification  \cite{K} of Laguerre planes.
\end{cor}
\begin{cor}\label{c3.3}
The group space $\mathbf{V}(\Delta(p,K)$ satisfies the condition $(\mathrm{Pgm})$.
\end{cor}
\begin{proof}
The group $\mathbf{T}(p,\Delta)$ is commutative (\cite{HP}, Theorem 4.14) and transitive. Hence the condition 5 is
satisfied (\cite{HP}, Proposition 6.5).
\end{proof}
\begin{cor}\label{c3.4}
Lines $A,B$ of the residual skewaffine plane are parallel iff
there exists a translation $\phi$ such that $\phi(A)=B$.
\end{cor}
\section{Some properties of residual skewaffine plane and their applications to Laguerre plane}
\begin{prop}\label{p4.1}
There are no three circles $L,M,N$ that are tangent at different
points with $L\in(p,K)$ and $M\cap N\subset \overline{p}$.
\end{prop}
\begin{proof}
Suppose circles $M,N$ are tangent to the circle $L$ of the pencil
$\langle p,K\rangle$ in points $x,y$ respectively and $M,N$ have
common point on the generator $\overline{p}$. By Corollary
\ref{c3.1} there exists $r\in L\setminus\overline{p}$ such that
$\mathrm{S}_{\overline{r},\overline{p};K}(x)=y$. We obtain
$\mathrm{S}_{\overline{r},\overline{p};K}(M)=N$ and hence $rM=rN$
is another common point of $M,N$.
\end{proof}
\begin{cor}\label{c4.1}
Parallel lines of $\mathbb{SA}(p,K)$ determined by circles with base points on a straight line determined by a circle
have a common point.
\end{cor}
\begin{prop}\label{p4.2}
Proper parallel lines of $\mathbb{SA}(p,K)$ determined by circles are disjoint iff their base points are distinct and
parallel.
\end{prop}
\begin{proof}
$\Leftarrow$ Let $m,n$ be distinct and parallel base points of circles $M,N$ and $M,N$ have a common point on the
generator $\overline{p}$. A translation $\alpha\in\mathbf{T}(p,\mathcal{G})$ such that $\alpha(m)=n$ transforms $M$ on
$N$, hence
$(M\setminus\overline{p})\cap(N\setminus\overline{p})=\emptyset$.\\
$\Rightarrow$ Assume the circles $M,N$ are tangent in a point
$q\in\overline{p}$ and their base points $m,n$ are not parallel.
Denote $L=(p,K,m)^{\circ}$, $z=nL$ and $M'=(z,L,q)^{\circ}$. By
Corollary \ref{c4.1} $M'$ is not tangent to $M$. This is a
contradiction with the proved part of the proposition.
\end{proof}
\begin{prop}\label{p4.3}
Suppose that base points of parallel lines $A,A'$ belongs to a straight
line $B$ determined by a circle. If a straight line $C$ determined
by a circle intersects the line $A$ then it intersects the line
$A'$.
\end{prop}
\begin{proof}
Let $x,x'$ be base points of lines $A$ and $A'$, respectively.
Let also $y$ be one of the common points of the lines $A,C$.
There exists a
translation $\tau\in\mathbf{T}(p,K)$ such that $\tau(x)=x'$. We
obtain $\tau(A)=A'$ and the point $y'=\tau(y)$ is a common point
of the lines $A',C$.
\end{proof}
\begin{lem}\label{l4.1}
Let circles $P,Q,R$ are tangent to a circle $L\in\langle p,K\rangle$ in pairwiese distinct points. If $Q$ have common
points with $P,R$, then $P,R$ have a common point.
\end{lem}
\begin{proof}
If two of the circles $P,Q,R$ determine parallel lines, the assertion follows from Proposition \ref{p4.3}. Assume $x\in
P\cap Q$, $y\in Q\cap R$ and $x,y\notin\overline{p}$. If $Q\in\langle p,K\rangle$, consider a translation
$\tau\in\mathrm{T}(p,K)$ such that $\tau(x)=y$ and the circle $P'=\tau(P)$. We obtain $y\in P'\cap R$ and $P'$ is
tangent to $L$. Hence, by Veblen-condition (V),the circles $P, R$ have a common point. In the case $Q\cap
L=\{r\}\neq\{p\}$ instead the translation $\tau$ we use the $\mathbb{L}$-strain $\phi\in\Delta(\overline{p},K,r)$ such
that $\phi(x)=y$ and the assertion follows by Veblen condition or Proposition \ref{p4.3}.
\end{proof}
\begin{defn}\label{d4.1}
Let $L$ be a circle $L$ of the pencil $\langle p,K\rangle$.
We say that points $a,b\not\in L$ are {\it equivalent} (under $L$)
and write $a\equiv_{L}b$ if
$|P\cap Q|\geq 1$ for any circles $P,Q$ tangent to $L$
and passing through $a$ and $b$, respectively.
\end{defn}
Lemma \ref{l4.1} gives the following.
\begin{prop}\label{p4.4}
For any circle $L\in\langle p,K\rangle$ the relation $\equiv_{L}$ is an equivalence on the set $\mathcal{P}\setminus
L$.
\end{prop}
\begin{prop}\label{p4.5}
For any points $a,b\in\mathcal{P}\setminus L$, $a\equiv_{L}b$ iff
there exist circles $P,Q$ tangent to $L$ in distinct point and
passing through $a,b$ respectively such that $|P\cap Q|\geq 1$.
\end{prop}
The points set of a special line can be described by the relation
$\equiv$ as follows:
\begin{prop}\label{p4.6}
Let $x\sqcup y$ be a special line, and let $L=(p,K,x)^{\circ}$, $z-y$, $z\neq x$. Then
$z\in x\sqcup y$ iff $z\equiv_{L}y$.
\end{prop}
\begin{proof}
$\Rightarrow$ Assume $y'\in(p,K,y)^{\circ}$, $M=(x,L,y')^{\circ}$
for some $y'\neq p,y$. By the definition of $x\sqcup y$ there
exists $\sigma\in\Delta(p,K)_{x}$ such that $\sigma(y)=z$. Then
$\sigma(M)=M$. If $z'=\sigma(y')$, we obtain
$y\equiv_{L}y'\equiv_{L}z'\equiv_{L}z$.\\
$\Leftarrow$ Let $y\in N\in\langle p,K\rangle$ and $M$ be a circle
passing through $z$ tangent to $L$ in a point different from $p$.
By $z\equiv_{L}y$ there exists a point $r$ such that $r\in M\cap
N$. Denote $P=(x,L,r)^{\circ}$ and $Q=(p,K,z)^{\circ}$. Then $P$
and $Q$ have a common point $s$ because $r\equiv_{L}z$. A strain
$\phi\in\Delta(\overline{p},K,x)$ such that $\phi(r)=s$ transforms
$y$ on $z$.
\end{proof}
\begin{prop}\label{p4.7}
If $x-p$, $x\neq p$, then for any point $y\notin L$, $x\equiv_{L}y$ iff there exist exactly two circles $M,M'$ tangent
to $L$ such that $x,y\in M,M'$.
\end{prop}
\begin{proof}
It is sufficient to prove $\Rightarrow$. Let $N$ be any circle tangent to $L$ such that $x\in N$ and
$P=(p,L,y)^{\circ}$. From $x\equiv_{L}y$ it follows that there exists $z\in P\cap N$. Then $M=\tau(N)$ where
$\tau\in\mathrm{T}(p,K)$, $\tau(z)=y$. Then we obtain $M'=S_{\overline{p},\overline{y},K}(M)$. Suppose, contrary to our
claim that there exists a circle $M''$ through $x,y$, tangent to $L$ and distinct from $M,M'$. Denote $r,r',r''$ the
base points of circles $M,M',M''$ respectively. There exists $\phi\in\Delta(p,K)_{r''}$ with $\phi(r)=r'$. We have:
$\phi(M'')=M''$, $\phi(M)=M'$ and $\phi(y)\neq x,y$. Hence $x,y,\phi(y)$ are three distinct points of two distinct
circles $M',M''$, a contradiction.
\end{proof}
\begin{lem}\label{l4.2}
For any $q$ parallel to $p$ and different from $p$ there exists exactly one $q'$ parallel to $p$ with the property:\\
$\forall x,y((x-\!\!\!\!\!/y,p-\!\!\!\!\!/x,y)\wedge
q\in(x(p,K,x)^{\circ},y)^{\circ}\rightarrow
q'\in(y,(p,K,y)^{\circ},x)^{\circ})$.
\end{lem}
\begin{proof}
The assertion is a consequence of the symmetry condition (P2). If the point $q$ determines the class of lines parallel
to a line $x\sqcup y$, then $q'$ determines the class of lines parallel to $y\sqcup x$.
\end{proof}
\begin{thm}\label{t4.1}
Let $q\neq p$, $q-p$, $x-\!\!\!\!\!/p$. Then the points of tangency of circles of the pencil $\langle p,K\rangle$ with
circles of the pencil $\langle x,q\rangle$ form a circle (without a point of the generator $\overline{p}$).
\end{thm}
\begin{proof}
Let $L=(p,K,x)^{\circ}$ and $M=(x,L,q)^{\circ}$. The point $x$ is the point of tangency of circles $L$ and $M$ of
pencils $\langle p,K\rangle$ and $\langle q,x\rangle$, respectively. Consider an arbitrary circle $N\in\langle
q,x\rangle$, $N\neq M$. By the condition (A3), there exists exactly one circle $P\in\langle p,K\rangle$ tangent to $N$
in some point $y$. The circle $Q=(x,L,y)^{\circ}$ is fixed by the group $\Delta(p,K)_{x}$. According to Corollary
\ref{c2.1}, any point of $Q$ distinct from $x$ and $pQ$ is an image of the point $y$ in some $\sigma\in\Delta(p,K)_{x}$.
Hence it is a point of tangency of circles of the pencils $\langle p,K\rangle$ and $\langle pQ,x\rangle$, respectively.
It follows that the circle $Q$ satisfies the assertion of the theorem.
\end{proof}
\begin{cor}\label{c4.2}
The circle $Q$ determined in Theorem \ref{t4.1} pass through the point $q'$ from Lemma \ref{l4.2}
\end{cor}
In the case of miquelian Laguerre planes of the characteristic different from 2,
the point $p$ in Theorem \ref{t4.1} can be chosen
arbitrarily by Remark \ref{r2.1}. For such planes we also obtain
the condition deciding whether a circle through two points tangent
to a circle can be constructed.
\begin{thm}\label{t4.2}
Let $x,y$ be points and $L$ a circle with $x-\!\!\!\!\!/y$, $x,y\notin L$ of miquelian Laguerre plane of characteristic
distinct from 2. The following conditions are equivalent:
\begin{enumerate}
\item[(1)] There exist exactly two circles through $x,y$ tangent to $L$.
\item[(2)] Any circle through $x$ tangent to $L$ intersects any circle through $y$ tangent to $L$.
\item[(3)] There exist two intersecting circles tangent to $L$ in distinct points containing $x,y$, respectively.
\end{enumerate}
\end{thm}
\begin{proof}
According to Remark \ref{r2.1} the assertion follows by Definition \ref{d4.1}, Proposition \ref{p4.5} and Proposition
\ref{p4.7} applied to the pencil $\langle xL,L\rangle$.
\end{proof}

\begin{rem}\label{r4.1}
In miquelian  Laguerre planes over a field $\mathbb{F}$ of
characteristic different from 2 the conditions of Theorem
\ref{t4.2} define the relation "$\equiv_{L}$" for any circle $L$.
In an analytic representation of such planes for a circle
$K=\{(x,0)|x\in\mathrm{F}\}\cup\{(\infty,0)\}$ points
$(a_{1},b_{1})$ and $(a_{2},b_{2})$ are equivalent with respect to
$K$ iff $b_{2}\in b_{1}\mathrm{F}^{2}$. In this case
classes of parallelity of special lines corresponds to classes of
squares of $\mathbb{F}$.
\end{rem}
\begin{rem}\label{r4.2}
If $\mathbb{F}$ is quadratically closed, then any special line
coincides with a generator and is a straight line. In this case
$\mathbb{SA}(p,K)$ contains two families of straight lines as
nearaffine residual plane connected with Minkowski plane (cf.
\cite{W}). But the class of straight lines determined by the circles of
the pencil $\langle p,K\rangle$ do not satisfy the condition about
existence of exactly one common point with other lines.
\end{rem}
% ----------------------------------------------------------------

%\bibliographystyle{amsplain}
%\bibliography{}
%%%%%%%%%%%%%%%%%%%%%%%%%%%%%%%%%%%%%%%%%%%%%%%%%%%%%%%%%%%%%%%%%%%%%%%%%%%%%%%%%%%%%%%%%%%%%%%%%%%%%%%%%%
\end{document}